\documentclass{amsart}
\usepackage{graphicx}
\usepackage{amssymb}
\usepackage{enumitem}

\newtheorem{thm}{Theorem}[section]

\theoremstyle{definition}

\theoremstyle{remark}
\newtheorem{rem}[thm]{Remark}

\begin{document}

\title[nowhere vanishing Hessian determinant]{Homogeneous functions with nowhere vanishing Hessian determinant}
\author{Connor Mooney}
\address{Department of Mathematics, UC Irvine}
\email{\tt mooneycr@math.uci.edu}

\begin{abstract}
We prove that functions that are homogeneous of degree $\alpha \in (0,\,1)$ on $\mathbb{R}^n$ and have nowhere vanishing Hessian determinant cannot change sign.
\end{abstract}
\maketitle

\section{Introduction}

Let $n \geq 2$, and let $\Omega \subset \mathbb{R}^n$ be the cone over a domain $\Sigma \subset \mathbb{S}^{n-1}$ that has nonempty boundary. Let 
$$p_n := \max\{1,\, (n-1)(n-2)\}.$$ 
In this paper we show:

\begin{thm}\label{Main}
If there exists a function $u : \Omega \rightarrow \mathbb{R}$ that satisfies:
\begin{enumerate}[label=(\roman*)]
\item $u$ is homogeneous of degree $\alpha \in (0,\,1)$,
\item $u \in W^{2,\,p}_{loc}(\Omega) \cap W^{1,\,\infty}((\Omega \cap B_2) \backslash B_1)$ for some $p > p_n$ (with $p = p_n$ allowed if $n \leq 3$),
\item $u > 0$ in $\Omega$ and $u = 0$ on $\partial \Omega$, and
\item either $\det D^2u$ or $-\det D^2u$ is locally strictly positive in $\Omega$,
\end{enumerate}
then $\mathbb{R}^n \backslash \Omega$ is a convex cone, and $\Sigma$ contains a closed hemisphere.
\end{thm}

\noindent By locally strictly positive we mean bounded below by positive constants almost everywhere on compact sets, where the constants may depend on the sets.

An immediate consequence of Theorem \ref{Main} is:

\begin{thm}\label{Main2}
If $u: \mathbb{R}^n \rightarrow \mathbb{R}$ satisfies:
\begin{enumerate}[label=(\roman*)]
\item $u$ is homogeneous of degree $\alpha \in (0,\,1)$,
\item $u \in W^{2,\,p}_{loc}\left(\mathbb{R}^n \backslash \{0\}\right)$ for some $p > p_n$, and
\item either $\det D^2u$ or $-\det D^2u$ is locally strictly positive in $\mathbb{R}^n \backslash \{0\}$,
\end{enumerate}
then $u$ does not change sign.
\end{thm}
\noindent Indeed, if $u$ changes sign then we may apply Theorem \ref{Main} to the sets $\{u > 0\}$ and $\{-u > 0\}$ to get a contradiction.

\begin{rem} 
Theorem \ref{Main2} is special to the cases $\alpha \in (0,\,1)$. Indeed, when $\alpha \notin [0,\,1]$ and $\alpha < k^2$ for some nonzero integer $k$, the $\alpha$-homogeneous functions
$$u = r^{\alpha}\cos(k\theta)$$
are sign-changing and have nowhere vanishing Hessian determinant on $\mathbb{R}^2 \backslash \{0\}$. We also remark that $0$-homogeneous functions have vanishing Hessian determinant on the rays where they achieve their maxima, and $1$-homogeneous functions have identically vanishing Hessian determinant. 
\end{rem}

Apart from its own interest, Theorem \ref{Main2} is motivated by the question of when interior gradient estimates hold for solutions to the special Lagrangian equation
\begin{equation}\label{sLag}
F(D^2u) := \sum_{k = 1}^n \tan^{-1}(\lambda_k(D^2u)) = \Theta(x) \in \left(-n\frac{\pi}{2},\, n\frac{\pi}{2}\right).
\end{equation}
Here $u$ is a function on a domain in $\mathbb{R}^n$ and $\lambda_k(D^2u)$ denote the eigenvalues of $D^2u$. Equation (\ref{sLag}) prescribes the mean curvature of the gradient graph of $u$ in $\mathbb{R}^n \times \mathbb{R}^n$. In particular, this graph is volume-minimizing when $\Theta$ is constant. 
The existence of continuous viscosity solutions to the Dirichlet problem for (\ref{sLag}) is known in certain situations (see e.g. \cite{HL}, \cite{CP}), and there are many fascinating open questions concerning the regularity of these solutions. For example, is not known whether they are locally Lipschitz if either $\Theta$ is a constant with $|\Theta| < (n-2)\frac{\pi}{2}$ (they are if $|\Theta| \geq (n-2)\frac{\pi}{2}$, see \cite{WY}) or if $\Theta$ is Lipschitz. Classical proofs of interior gradient estimates for elliptic PDEs involve differentiating the equation once, so is reasonable to ask if interior gradient estimates for (\ref{sLag}) hold under such conditions on $\Theta$.
A first attempt to {\it disprove} the validity of such estimates could be to build a function $u$ that is homogeneous of degree $\alpha \in (0,\,1)$, smooth away from the origin, and has nowhere vanishing Hessian determinant. Then $F(D^2u)$ would behave near the origin like a constant plus smooth function on the sphere times $|x|^{2-\alpha}$, which is $C^1$, while $u$ has unbounded gradient. Taking $u = |x|^{\alpha}$ appears to do the trick, but this function is not a viscosity solution at the origin; one needs $u$ to change sign to prevent this issue. Theorem \ref{Main2} precludes the existence of such functions, and hence can be viewed as evidence in favor of a positive result.

\begin{rem}
One might also try to build one-homogeneous functions $u$ on $\mathbb{R}^n$ such that $\sigma_{n-1}(D^2u)$ is nowhere vanishing, since in that case $F(D^2u)$ is Lipschitz at the origin. Here $\sigma_k(D^2u)$ denotes the $k^{th}$ symmetric polynomial of the eigenvalues. It is not hard to show that such functions are necessarily convex or concave (see Section \ref{Preliminaries}), and thus do not solve the equation at the origin. Interestingly, there exist nonlinear one-homogeneous functions $u$ on $\mathbb{R}^3$ whose Hessians are either indefinite or $0$ at every point (see \cite{M}), so that $F(D^2u)$ tends to zero at the origin along rays, but $F(D^2u)$ is not continuous at the origin for these examples.
\end{rem}

The paper is organized as follows. In Section \ref{Preliminaries} we recall some preliminary results about one-homogeneous functions and about maps with integrable dilatation, which are natural analogues of quasi-conformal maps in higher dimensions. In Section \ref{Proof} we prove Theorem \ref{Main}. The idea of the proof is to study the geometry of the gradient image of the one-homogeneous function $u^{\frac{1}{\alpha}}$. Our analysis is partly inspired by the beautiful arguments in \cite{HNY} used to show the linearity of one-homogeneous functions on $\mathbb{R}^3$ that solve linear uniformly elliptic equations.

\section*{Acknowledgments}

This research was supported by NSF grant DMS-1854788.

\section{Preliminaries}\label{Preliminaries}

In this section we recall a few results about one-homogeneous functions and about maps with integrable dilatation.

\vspace{2mm}

Let $v$ be a one-homogeneous function on $\mathbb{R}^n$ that, away from the origin, is locally $W^{2,\,p}$ for some $p \geq 1$. Euler's formula for homogeneous functions says that
\begin{equation}\label{Euler}
v(x) = \nabla v(x) \cdot x.
\end{equation}
Here and below we let $r := |x|$ and we denote points in $\mathbb{S}^{n-1}$ by $\omega$. Writing 
$$v = r\,g(\omega)$$
and choosing a coordinate system where $\omega$ is the last direction, we have
\begin{equation}\label{OneHomogHessian}
D^2v(\omega) = 
 \begin{bmatrix}
     \nabla_{\mathbb{S}^{n-1}}^2g + g\,I_{n-1 \times n-1} & 0 \\ 
     0 & 0
 \end{bmatrix}.
\end{equation}
Here and below, $\nabla_{\mathbb{S}^{n-1}}$ and $\nabla_{\mathbb{S}^{n-1}}^2$ denote the usual gradient and Hessian operators on the sphere. It is sometimes useful to represent $v$ in $\{x_n > 0\}$ by a function $\bar{v}$ on $\mathbb{R}^{n-1}$ defined by
$$\bar{v}(y) := v(y,\,1),$$
so that
$$v(x',\,x_n) = x_n\bar{v}\left(\frac{x'}{x_n}\right).$$
Taking the Hessian yields 
$$D^2v(y,\,1) = 
 \begin{bmatrix}
     D^2\bar{v} & -D^2\bar{v} \cdot y \\ 
     -D^2\bar{v} \cdot y & y^T\cdot D^2\bar{v} \cdot y
 \end{bmatrix}.
 $$
Using this we calculate
\begin{equation}\label{DetRelation}
\sigma_{n-1}(D^2v)(y,\,1) = \text{tr}(\text{cof}(D^2v))(y,\,1) = (1+|y|^2)\det D^2\bar{v}(y).
\end{equation}
Here the operator $\sigma_k$ denotes the $k^{th}$ symmetric polynomial of the eigenvalues. It is easiest to verify this formula after rotating in the $y$ variables so that $D^2\bar{v}$ is diagonal.
If $v$ is $C^2$ in a neighborhood of $e_n$ and $D^2\bar{v}(0)$ is nonsingular, then we can represent the gradient image of $v$ near $e_n$ as the graph of a function $w$ using the relation
$$w(\nabla \bar{v}(y)) = \partial_nv(y,\,1) = \bar{v} - y \cdot \nabla \bar{v},$$
i.e. $w$ is the (negative) Legendre transform of $\bar{v}$. One differentiation gives
$$\nabla w(\nabla\bar{v}(y)) = -y,$$
and another gives
\begin{equation}\label{II}
D^2w(\nabla \bar{v}(y)) = -(D^2\bar{v})^{-1}(y).
\end{equation}
In particular, the second fundamental form at $\nabla v(e_n)$ of the image under $\nabla v$ of a small ball around $e_n$ is $(D^2\bar{v})^{-1}(0)$. From this it is easy to see that if $v$ is locally $C^2$ away from the origin and $\sigma_{n-1}(D^2v)$ is nowhere vanishing, then $v$ is either convex or concave. Indeed, it suffices to show that either $D^2v \geq 0$ or $-D^2v \geq 0$ at some point. By (\ref{II}) this is true at the inverse image under $\nabla v$ of any point on $\nabla v(\mathbb{S}^{n-1})$ that is touched from one side by a hyperplane.

\vspace{2mm}

We now recall a few facts about maps of integrable dilatation. Let 
$$\varphi: U \subset \mathbb{R}^n \rightarrow \mathbb{R}^n$$ 
be a map in $W^{1,\,n}_{\text{loc}}(U)$ such that $\det D\varphi > 0$ almost everywhere. The dilatation $K$ of $\varphi$ is defined by the ratio
$$K(x) := \frac{|D\varphi|^n}{\det D\varphi}.$$
If $K$ is bounded and $n = 2$ then $\varphi$ is quasi-conformal, hence continuous and either open or constant by classical results. Reshetnyak extended this result to mappings with bounded dilatation in all dimensions \cite{R}. The boundedness required in Reshetnyak's theorem has since been relaxed to integrability in certain $L^p$ spaces:

\begin{thm}[{\bf Iwaniec-\v{S}ver\'{a}k}, \cite{IS}]\label{OM1}
If $n = 2$ and $K \in L^1_{loc}(U)$, then $\varphi$ is continuous and either open or constant.
\end{thm}

\begin{thm}[{\bf Manfredi-Villamor}, \cite{MV}]\label{OM2}
If $n \geq 3$ and $K \in L^p_{loc}(U)$ for some $p > n-1$, then $\varphi$ is continuous and either open or constant.
\end{thm}

\noindent It is conjectured that the latter result holds in the equality case $p = n-1$ (see \cite{IS}), and there are counterexamples when $p < n-1$ due to Ball (see \cite{B}).

\section{Proof of Theorem \ref{Main}}\label{Proof}

In this final section we prove the main theorem.

\begin{proof}[{\bf Proof of Theorem \ref{Main}}]
We write
$$u = r^{\alpha}f(\omega).$$
In a coordinate system where the last direction is $\omega \in \Sigma$, the Hessian of $u$ at $\omega$ can be written
\begin{equation}\label{alphaHessian}
D^2u = \begin{bmatrix}
     \nabla_{\mathbb{S}^{n-1}}^2f + \alpha\,f\,I_{n-1 \times n-1} & (\alpha - 1)\nabla_{\mathbb{S}^{n-1}}f \\ 
     (\alpha - 1)\nabla_{\mathbb{S}^{n-1}}f & \alpha(\alpha - 1)f
 \end{bmatrix}.
 \end{equation} 
Subtracting the multiple $\frac{\left(\nabla_{\mathbb{S}^{n-1}}f\right)_k}{\alpha f}$ of the last row from the $k^{th}$ row in (\ref{alphaHessian}) for $k \leq n-1$ and taking the determinant we arrive at
\begin{equation}\label{uHessDet}
\begin{split}
\det D^2u &= \alpha(\alpha -1)f^{2-n} \cdot \\
&\left[\det\left(f\nabla_{\mathbb{S}^{n-1}}^2f + \left(\frac{1}{\alpha}-1\right)\nabla_{\mathbb{S}^{n-1}}f \otimes \nabla_{\mathbb{S}^{n-1}}f + \alpha f^2I_{n-1 \times n-1}\right)\right].
\end{split}
\end{equation}
Now let $v$ be the one-homogeneous function defined by
$$v = \begin{cases}
u^{1/\alpha}, \quad \omega \in \Sigma \\
0, \quad \text{otherwise.}
\end{cases}$$
At $\omega \in \Sigma$ we compute the Hessian of $v$ in the same coordinates as above, using the formula (\ref{OneHomogHessian}):
 \begin{equation}\label{vHessDet0}
 D^2v = \frac{1}{\alpha}f^{\frac{1}{\alpha}-2}
 \begin{bmatrix}
     f\nabla_{\mathbb{S}^{n-1}}^2f + \left(\frac{1}{\alpha}-1\right)\nabla_{\mathbb{S}^{n-1}}f \otimes \nabla_{\mathbb{S}^{n-1}}f + \alpha f^2I_{n-1 \times n-1} & 0 \\ 
     0 & 0
 \end{bmatrix}.
 \end{equation}
We conclude from (\ref{uHessDet}) and (\ref{vHessDet0}) that
\begin{equation}\label{vHessDet}
 \sigma_{n-1}(D^2v) = \frac{f^{\frac{n-1}{\alpha} - n}}{\alpha^n(\alpha-1)}\det D^2u
\end{equation}
on $\Sigma$. In particular, either $\sigma_{n-1}(D^2v)$ or $-\sigma_{n-1}(D^2v)$ is locally strictly positive in $\Omega$. We also have by standard embeddings and the fact that 
$$\frac{1}{\alpha} > 1$$ 
that $v \in C^1(\mathbb{S}^{n-1})$. Indeed, by homogeneity we may view $v$ as a function of $n-1$ variables away from $0$. When $n \geq 4$ the Sobolev exponent $p_n$ is thus supercritical. In the case $n = 3$ it is critical and, denoting by $\bar{v}$ the restriction of $v$ to a hyperplane tangent to $\mathbb{S}^{n-1}$, we may apply the continuity assertion in Theorem \ref{OM1} to either $\nabla \bar{v}$ or its reflection over a line. Here we used the fact that $\sigma_{n-1}(D^2v)$ is nowhere vanishing in $\{v > 0\}$ and the relation (\ref{DetRelation}). In the case $n = 2$ we use that $W^{1,\,1}$ embeds to continuous on the line.
 
Now let $K := \nabla v\left(\mathbb{S}^{n-1}\right)$. For $\nu \in \mathbb{S}^{n-1}$, slide the hyperplane $\{x \cdot \nu = t\}$ (starting with $t$ large, and decreasing $t$) until it touches $K$ at some point $p_{\nu}$. Since $0 \in K$ we have $p_{\nu} \cdot \nu \geq 0.$ We claim that
\begin{equation}\label{PositiveSet}
p_{\nu} \cdot \nu > 0 \Rightarrow (\nabla v)^{-1}(p_{\nu}) \cap \mathbb{S}^{n-1} = \{\nu\}.
\end{equation}
To show the implication (\ref{PositiveSet}) it suffices to show that $(\nabla v)^{-1}(p_{\nu}) \cap \mathbb{S}^{n-1} \subset \{\nu,\,-\nu\}$, since by (\ref{Euler}) we have 
$$\nabla v(\omega) \cdot \omega = v(\omega) \geq 0$$ 
for all $\omega \in \mathbb{S}^{n-1}$. 
Assume by way of contradiction that $p_{\nu} \cdot \nu > 0$ but $(\nabla v)^{-1}(p_{\nu}) \cap \mathbb{S}^{n-1}$ is not contained in $\{\nu,\,-\nu\}$. After a rotation we may assume that $\nu = e_1$, and after another rotation in the $x_2,\,...,\,x_n$ variables we may assume that $\nabla v(\tilde{\nu}) = p_{\nu}$ for some $\tilde{\nu} \in \mathbb{S}^{n-1}$ such that $\tilde{\nu}_n > 0$. For $y \in \mathbb{R}^{n-1}$ we let
$$\bar{v}(y) = v(y,\,1).$$
By construction, $\partial_1\bar{v}$ has a local maximum at $\tilde{\nu}/\tilde{\nu}_n$ (here we identify points on the hyperplane $\{x_n = 1\}$ with $\mathbb{R}^{n-1}$), and $\bar{v} > 0$ at this point, since $\nabla v(\tilde{\nu}) = p_{\nu}$ is nonzero. However, using (\ref{DetRelation}) we see that
either $\det D^2\bar{v}$ or $-\det D^2\bar{v}$ is locally strictly positive in $\{\bar{v} > 0\}$. Theorems (\ref{OM1}) and (\ref{OM2}), applied to either $\nabla \bar{v}$ or its reflection over a hyperplane, imply that $\nabla \bar{v}$ is an open mapping in $\{\bar{v} > 0\}$, which contradicts that $\partial_1\bar{v}$ has a local maximum in this set.

Finally, let $\text{co}(K)$ denote the convex hull of $K$, and let $w$ be the support function of $\text{co}(K)$, that is,
$$w(x) := \sup_{y \in \text{co}(K)} (y \cdot x).$$ 
The implication (\ref{PositiveSet}) implies that $v = w$. Indeed, it is clear that $0 \leq v \leq w$, and for $\nu \in \mathbb{S}^{n-1} \cap \{w > 0\}$ we have
$$w(\nu) = p_{\nu} \cdot \nu = \nabla v(\nu) \cdot \nu = v(\nu).$$ 
Because either $\sigma_{n-1}(D^2v)$ or $-\sigma_{n-1}(D^2v)$ is locally strictly positive in $\Omega$, the set $\text{co}(K)$ has non-empty interior. Indeed, if not, then $v = w$ is translation-invariant in some direction orthogonal to $\text{co}(K)$, which along with the one-homogeneity of $v$ implies that $\sigma_{n-1}(D^2v) \equiv 0$. We conclude that $\{v > 0\} \cap \mathbb{S}^{n-1}$ contains some closed hemisphere, completing the proof.
\end{proof}

\begin{figure}
 \begin{center}
    \includegraphics[scale=0.5, trim={0mm 65mm 0mm 20mm}, clip]{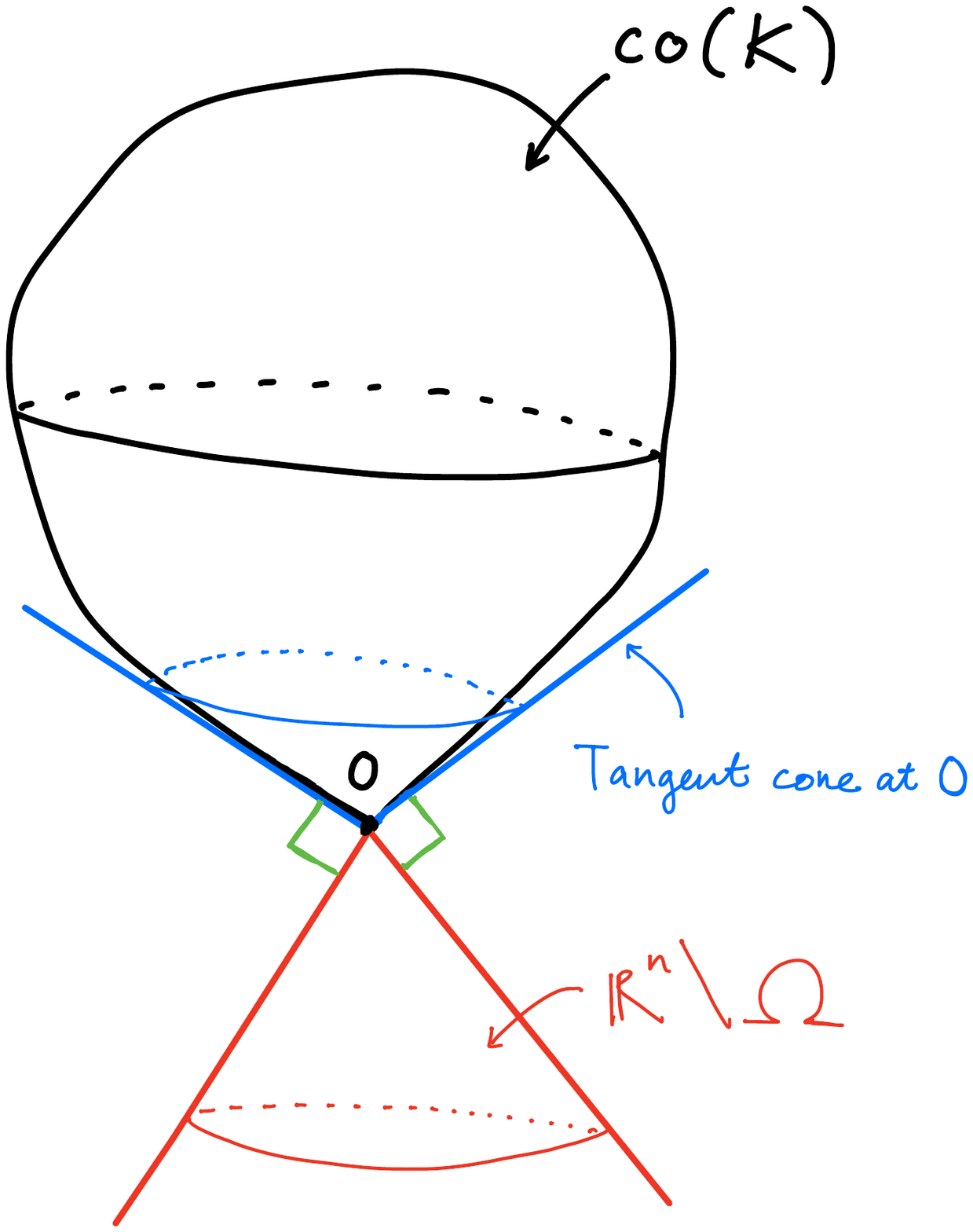}
\caption{}
\label{fig1}
\end{center}
\end{figure}

\begin{rem}
The proof in fact shows that $u$ is the $\alpha^{th}$ power of the support function of a convex set that has nonempty interior and $0$ in its boundary (note that if $0$ were not in the boundary of $\text{co}(K)$ then $w = v > 0$ on all of $\mathbb{S}^{n-1}$), and that $\mathbb{R}^n \backslash \Omega$ is the reflection through the origin of the convex dual to the tangent cone of this set at the origin (see Figure \ref{fig1}).
\end{rem}



\end{document}